\input amstex

\documentstyle{amsppt}

\NoBlackBoxes

\magnification 1200

\hsize 6.5truein

\vsize 8truein

\topmatter

\title The Grone Merris conjecture and a quadratic eigenvalue problem \endtitle

\author Nets Hawk Katz \endauthor

\affil Indiana University \endaffil
\thanks  
 The author was partially supported by a National Science Foundation
Grant \endthanks

\endtopmatter

\head \S 0 Introduction \endhead

Given a graph $G=(V,E)$, we define the transpose degree sequence $d_j^T$ to
be equal to the number of vertices of degree at least $j$. We define $L_G$, the
graph Laplacian, to be the matrix, whose rows and columns are indexed by the vertex
set $V$, whose diagonal entry at $v$ is the degree of $v$ and whose value at a pair $(v,w)$ is $-1$ if 
$(v,w) \in E$ and $0$ otherwise.

Grone and Merris conjectured [GM]

\proclaim{Conjecture}
If $\lambda_1,\dots,\lambda_m$ are the eigenvalues of $L_G$ in (weakly) decreasing order,
then for any $1 \leq j \leq m$, we have
$$\sum_{l=1}^j \lambda_l \leq \sum_{l=1}^j d_l^T.$$
\endproclaim.

The first inequality is well known and the last inequality is indeed always an equality. On the
class of threshold graphs, all the inequalities are equalities. The second inequality was proved by
Duval and Reiner [DR], and this paper grew out of an attempt to understand their proof.

We will say that a graph is semi-bipartite if its vertex set is the union of a clique and an isolated
set. We say a semi-bipartite graph is $k$ regular, if every vertex in the isolated set has degree
$k$. Duval and Reiner proved the second inequality by observing that it suffices to prove it
for $1$-regular semi-bipartite graph. They proved the second inequality by showing it was trivial
for all but a few classes of $1$ regular semi-bipartite graphs and then solving those cases one
by one using a computer algebra system. However, if $j$ is the number of vertices in the clique
adjacent to some vertex in the isolated set, then it seems that the $j$th inequality is quite
difficult to prove and often fails to be an equality very narrowly. The $j$th inequality rather 
easily implies all the others.

The purpose of this paper is to prove

\proclaim{Main Theorem} The Grone-Merris conjecture holds for 1-regular semibipartite graphs.
\endproclaim

The prood consists of an analysis of the roots of a certain polynomial. This polynomial
arises from a quadratic eigenvalue problem (see [TM]  for more information on QEP's)
which is equivalent to finding the spectrum of the Laplacian. The analysis has two
parts. The first part involves writing down lower degree polynomials each of which vanishes
at a root of the original polynomial and to show that the lower degree polynomials are totally
ordered by majorization. The second part is to use a homotopy to relate the order of root
in the lower degree polynomials to the order of roots in the original polynomial.

We believe there is a good chance of extending the results of this paper to the setting of general
semi-bipartite graphs. Every semi-bipartite graph gives rise to a polynomial eigenvalue problem
in determining its spectrum. The associated lower degree polynomials are indeed totally
ordered with respect to majorization. However some technical details may need to be worked out
to ensure that the homotopy method works.

We add that semi-bipartite graphs are an extremely natural setting in which to study
the Grone Merris conjecture. In particular, all threshold graphs are semi-bipartite.Also
we have some intuition that edges towards higher degree vertices have more impact
on the spectrum of the Laplacian than edges connecting low degree vertices. Therefore,
we hope we are on a path leading towards resolution of the conjecture and we hope
this paper gives rise to future work.

The paper is organized as follows: In section1, we state preliminary lemmas in linear algebra.
In section 2, we derive our quadratic eigenvalue problem. In section 3, we give a couple of elementary
examples which demonstrate the difficulty of the inequalities and illustrate our method. In section 4,
we complete the analysis of the polynomial. In section 5, we do some bookkeeping which proves the
main theorem.

The author would like to thank Vic Reiner for telling him the problem and would like to thank
Hari Bercovici and  Elizabeth Housworth for helpful discussions.

\head \S 1 Preliminaries \endhead

Throughout this paper, we will use the conventional notation that if
$\lambda=(\lambda_1,\dots,\lambda_j)$ and $\eta=(\eta_1,\dots,\eta_j)$ are
(weakly) decreasing
vectors of real numbers then
$$\lambda \triangleright \eta,$$
whenever for any $1 \leq l \leq j$, we have that
$$\sum_{k=1}^l \lambda_k \geq \sum_{k=1}^l \eta_k.$$
We read this $\lambda$ majorizes $\eta$.

We will abuse this notation as follows. If $p$ and $q$ are polynomials of the same
degree with all real roots, we write $p \triangleright q$ provided that the roots
of $p$, written in decreasing order, majorize the roots of $q$.

We will abuse the notation further. If $A$ and $B$ are matrices of the same
dimension
with all eigenvalues real, we will write $A \triangleright B$ provided that
the eigenvalues of $A$ written in decreasing order majorize the eigenvalues of
$B$.

We recall a basic fact from linear algebra.

\proclaim{Proposition 1.1} Let $A$ be a symmetric $j \times j$ matrix. Let
$\lambda_1,\dots,\lambda_j$ be the eigenvalues in decreasing order. Then
$$\sum_{k=1}^l \lambda_k \geq \operatorname{trace}(EAE),$$
where $E$ is any orthogonal projection of rank $l$. The inequality is
attained when $E$ is the projection into the first $l$ eigenvectors.
\endproclaim

This has many important consequences. We state a few.

\proclaim{Proposition 1.2} Let $\lambda_1,\dots,\lambda_j$ be 
a decreasing sequence of real numbers.
Let $A$ be a diagonal matrix with the $\lambda$'s along the diagonal.
$$A=\pmatrix \lambda_1 & 0 & \dots & 0 \\
             0 & \lambda_2 & \dots & 0 \\
             \vdots & \vdots & \vdots & \vdots \\
             0 & 0 & \dots & \lambda_j \endpmatrix.$$
Let $B$ be a martix with zeroes along the diagonal and 1's everywhere else:
$$B=\pmatrix 0 & 1 & \dots & 1 \\
             1 & 0 & \dots & 1 \\
             \vdots & \vdots & \vdots &\vdots \\
             1 & 1 & \dots & 0 \endpmatrix.$$
Then 
$$A+B \triangleright A,$$
and
$$A-B \triangleright A.$$
\endproclaim

\demo{Proof} Let $E_l$ be the orthogonal projection of a vector into it first
$l$ components.
$$\lambda_1+\lambda_2 + \dots \lambda_l = \operatorname{trace}( E_l(A+B) E_l)
= \operatorname{trace}( E_l(A-B) E_l).$$
Now we apply Proposition 1.1 \qed \enddemo

\proclaim{Corollary 1.3} Let $a,b$ be positive real numbers with $b<a$. Then
with $A$ and $B$ as above
$$A-aB \triangleright A-bB.$$
(Similarly, we have
$$A+aB \triangleright A+bB.)$$
\endproclaim

\demo{Proof} Let $E$ be the $l$ dimensional orthogonal projection into the $l$
greatest eigenvalues of $A-bB$.
Then
$$\operatorname{trace}( E(A-bB)E)=\operatorname{trace}(EAE) -b
\operatorname{trace}(EBE).$$
Applying the argument of Proposition 1.2, we see that $A-bB$ majorizes $A$ so
that since $E$ achieves the maximum, it must be that $\operatorname{trace}(EBE)$ is
nonpositive. But 
$$\operatorname{trace}( E(A-aB)E)=\operatorname{trace}(EAE) -a
\operatorname{trace}(EBE)
\geq \operatorname{trace}( E(A-bB)E).$$
Proposition 1.1 yields the desired result. \qed \enddemo

We will meet the matrix $A+B$ frequently in this paper. In the following
proposition,
we compute its determinant.

\proclaim{Proposition 1.4} With $A$ and $B$ as above, and none of the
$\lambda_j$ equal to
1,
$$\det(A+B)=(\prod_{l=1}^j (\lambda_l-1)) (1+\sum_{l=1}^j {1 \over \lambda_l-1}).$$
In the case that one or more of the $\lambda_j$'s is equal to 1, the determinant
may be read from the above by formal cancellation. If two or more of the
$\lambda_j$'s equal 1, the determinant is 0.
\endproclaim

\demo{Proof}
$$A+B=A-I + (B+I).$$
Now just observe that $B+I$ has rank 1. The determinant is a multilinear
antisymmetric
function on the rows. Thus when we expand $\det(A-I + (B+I))$, any term
involving more than one row of $B+I$ disappears. We are left with the expression
above.
\qed \enddemo

\head \S 2 Graph Laplacians and quadratic eigenvalue problems \endhead

Let $G=(V,E)$ be graph. For each vertex $v \in V$, let $d_v$ be
its degree.

We define its Laplacian $L_G$, to be a matrix whose rows and columns
are indexed by $v$, and whose components are given by $(L_G)_{vv}=d_v$
while for $v,w \in V$ with $v \neq w$, we have
$$(L_G)_{vw} = -chi_E(v,w),$$
where $\chi_E$ is the indicator function of $E$, the set of edges.
The Laplacian $L_G$ is always positive definite since it is the matrix
associated to the quadratic form
$$Q_G(x,x) = \sum_{(v,w) \in E} (x_v-x_w)^2.$$

We restrict to a special class of graphs. We say that a graph 
$G=(V,E)$ is semi-bipartite if $V$ can be written as a disjoint
union $V=V_g \cup V_b$ where $V_g$ is a clique 
complete while $V_b$ is an isolated set. We say a semi-bipartite
graph $G=(V_g \cup V_b,E)$
is $d$-regular for a fixed positive integer $d$ provided
every vertex of $V_b$ has degree $d$. From this point on,
we will be considering a $1$-regular semibipartite graph.

Following Duval and Reiner, we divide our vertices into various
classes, permutations of which are automorphisms of the graph.
We denote  $\#(V_g)=n$. We denote by $j$, the number of vertices
of $V_g$ which are adjacent to some vertex of $V_b$. Clearly
$j \leq n$. We consider the case when $j=n$ as degenerate, and we 
defer its consideration to \S 5. For now, we assume $j<n$. 
We denote by $v_1, \dots, v_j$ those vertices of $V_g$ which are
adjacent to some vertex of $V_b$.
We let $V_{g,extra}$ to be the set of those $n-j$ vertices of $V_g$
not adjacent to any vertex in $V_b$. We denote by $W_l$, the set of vertices of $V_b$ which are adjacent
to the vertex $v_l$ . We denote by $k_l$ the
cardinality of $W_l$. 

We view $L_G$ as acting on functions on $V$, and we are now in a position
to identify several invariant subspaces of $L_G$. We say a function $f$
has sum zero if $\sum_{v \in V} f(v)=0$.

On any function supported on $V_{g,extra}$ having sum zero, the Laplacian
$L_G$ acts by multiplication by $n$. We have identified $n-j-1$ eigenvalues
equal to $n$. On any function supported on some $W_l$ 
which has sum zero, the Laplacian $L_G$ acts by multiplication by $1$.
We have identified $\#(V_b) - j$ eigenvalues equal to $1$.
To understand the spectrum of $L_G$, it now suffices to restrict to
the space of functions which are constant on $V_{g,extra}$ and on each 
$W_l$. We describe an orthonormal basis for this space.

Here for $1<l<j$,the vector $e_l$ will be the function equal to 1 on $v_l$ and zero elsewhere, $e_{extra}$ the vector
equal to ${1 \over \sqrt{n-j}}$ on $V_{g,extra}$ and zero elsewhere and
$f_l$, the vector equal to ${1 \over \sqrt{k_{l}}}$ on $C$ and zero elsewhere. We proceed to write down the $2j+1$
dimensional matrix of $L_G$ acting on this basis, which we will denote
(considering $j$ to be fixed) as $M(n,k_1,\dots,k_j)$.

$$M(n,k_1,\dots,k_j)=$$
$$\pmatrix n + k_1 - 1 & -1 &\dots & -1 & -\sqrt{n-j} & -\sqrt{k_1} & 0 & \dots & 0 \\
         -1 & n+k_2-1 & \dots & -1 &-\sqrt{n-j} & 0 & -\sqrt{k_2}
 & \dots & 0\\
          \vdots & \vdots & \vdots & \vdots & \vdots & \vdots
          &\vdots &\vdots &\vdots \\
         -1 & -1 & \dots & n+k_j-1 & -\sqrt{n-j} & 0 & 0 &\dots
 &-\sqrt{k_j} \\
          -\sqrt{n-j} &-\sqrt{n-j} &\dots &-\sqrt{n-j} &j & 0 & 0
&\dots & 0 \\
          -\sqrt{k_1} & 0 &\dots & 0 & 0 & 1 & 0 &\dots &0 \\
          0 & -\sqrt{k_2} & \dots & 0 & 0 & 0 & 1 &\dots & 0 \\
         \vdots &\vdots &\vdots &\vdots &\vdots &\vdots
&\vdots &\vdots &\vdots \\
         0 &0 &\dots &-\sqrt{k_j} &0 &0 &0 &\dots & 1 \endpmatrix $$

Duval and Reiner studied the cases $j=2$ and $j=3$ using computer
algebra. The remaining cases they avoided because they were only trying to prove the second inequality. We shall analyze this matrix directly.

First note that $M(n,k_1,\dots,k_j)$ has one eigenvalue which is zero
since the vector 
$$(1,\dots, 1, \sqrt{n-j},\sqrt{k_1},\dots,\sqrt{k_j})$$
is in the kernel. We reduce the dimension by one writing a basis for
the orthonormal complement to the kernel.

$$g_l=e_l -{1 \over \sqrt{n-j}} e_{extra},$$
and
$$h_l=f_l-{\sqrt{k_l} \over \sqrt{n-j}} e_{extra}.$$

In this basis, the Laplacian $L_G$ restricts to the matrix
$$N(n,k_1,\dots,k_j)=$$
$$\pmatrix n + k_1 & 0 &\dots & 0  & 0 & \sqrt{k_1} & \dots & \sqrt{k_1} \\
         0 & n+k_2 & \dots &  & \sqrt{k_2} & 0
 & \dots & \sqrt{k_2}\\
          \vdots & \vdots & \vdots  & \vdots & \vdots
          &\vdots &\vdots &\vdots \\
         0 & 0 & \dots & n+k_j  & \sqrt{k_j} & \sqrt{k_j} &\dots
 &0 \\
         -\sqrt{k_1} & 0 &\dots & 0  & 1 & 0 &\dots &0 \\
          0 & -\sqrt{k_2} & \dots & 0  & 0 & 1 &\dots & 0 \\
         \vdots &\vdots &\vdots  &\vdots &\vdots
&\vdots &\vdots &\vdots \\
         0 &0 &\dots &-\sqrt{k_j} &0  &0 &\dots & 1 \endpmatrix $$

We take a less magnified view of this matrix:
$$N(n,k_1,\dots,k_j)=\pmatrix A & B \\
                                                     C & I \endpmatrix.$$
 Here $A$,$B$, and $C$ are $j \times j$ matrices and $I$ is the identity. Now if $\lambda$ is an
 eigenvalue different from 1 of $N(n,k_1,\dots,k_j)$, it must be we can find a pair $v,w$ of vectors,
 at least one of which is nonzero so that
 $$(A-\lambda) v + B w= 0,$$
 and
 $$C v +(I-\lambda) w=0.$$
 Solving for $w$ in terms of $v$ and observing that $I-\lambda$ commutes with everything,
 we obtain that the matrix
 $$(A-\lambda)(I-\lambda)-BC$$
 is non-invertible. This is the promised quadratic eigenvalue problem. In our case, we see that
 $$(A-\lambda)(I-\lambda)-BC=$$
 $$\pmatrix   (n+k_1-\lambda)(1-\lambda) & k_2 & \dots & k_j \\
                        k_1 &  (n+k_2-\lambda)(1-\lambda) & \dots & k_j \\
                        \vdots & \vdots & \vdots &\vdots \\
                        k_1 & k_2 & \dots & (n+k_j-\lambda)(1-\lambda) \endpmatrix$$
 Dividing on the right by the matrix
 $$\pmatrix   k_1 & 0 & \dots & 0 \\
                       0 & k_2 & \dots & 0 \\
                       \vdots &\vdots &\vdots &\vdots \\
                       0 & 0 &\dots &k_j \endpmatrix,$$
 we see we must find the values where
 $$\pmatrix   {(n+k_1-\lambda)(1-\lambda) -k_1\over k_1} & 0 & \dots & 0 \\
                        0 &  {(n+k_2-\lambda)(1-\lambda) -k_2 \over k_2}& \dots & 0 \\
                        \vdots & \vdots & \vdots &\vdots \\
                        0 & 0 & \dots & {(n+k_j-\lambda)(1-\lambda) -k_j \over \lambda}\endpmatrix
   +\pmatrix  1 & 1 & \dots & 1\\
                     1 & 1 & \dots & 1\\
                     1 & 1 & \dots & 1\\   
                     1 & 1 & \dots & 1 \endpmatrix.$$  
  We can take the determinant of the above sum using Proposition 1.4. We obtain the polynomial
  $$F_{n,k_1,\dots,k_j}(\lambda)=G_{n,k_1, \dots k_j}(\lambda) \prod_{l=1}^{j}
[(n+k_l-\lambda)(1-\lambda)-k_l],$$
where
 $$G_{n,k_1, \dots k_j}(\lambda)=1 + \sum_{l=1}^j {k_l \over
(n+k_l-\lambda)(1-\lambda)-k_l}.$$

The sum of the $j$ largest roots of  $F_{n,k_1,\dots,k_j}(\lambda)$ shall be the subject of
sections 3 and 4.

\head \S 3 Illustrative Examples \endhead

Our goal in this section and the next is to bound the $j$ largest roots of $F_{n,k_1,\dots,k_j}(\lambda)$ 
by
$$jn +\sum_{l=1}^j k_l.$$
This bound is the $j$th inequality in the Grone Merris conjecture.

In this section, we do a couple of simple examples. The first is intended to demonstrate that our
inequality can seem very close to being tight even far from the threshold situation. The second is
intended to illustrate the method we shall develop in \S 4.

Our first example is the case when the $k$'s are all equal. That is
$$k_1=k_2 = \dots =k_j =k.$$
In that case, we can readily factor
$$F_{n,k_1,\dots,k_j}(\lambda)=( (n+k-\lambda)(1-\lambda)-k)^{j-1}     ((n+k-\lambda)(1-\lambda) +(j-1)k).$$

Here, we can write out the $j$ largest roots explicitly. We have a root with multiplicity $j-1$
given by
$$r_1={n + k + 1 + \sqrt{(n+k+1)^2-4n} \over 2},$$
and a root with multiplicity 1 given by
$$r_2={n + k + 1 + \sqrt{(n+k+1)^2-4n-4jk} \over 2}.$$
Observing that the square root function is concave, we see that
$$(j-1) r_1 + r_2 \leq j({n+k+1 + \sqrt{(n+k+1)^2 -4(n+k)} \over 2} = j (n+k),$$
which is exactly the desired inequality.

This shows that the inequality can be very tight: concavity of square roots is a second order effect.
It also suggests that the Grone Merris conjecture might be viewed as a type of convexity inequality,
although we have been unable to carry this out.

Our second example is one of those studied by Duval and Reiner using computer algebra.
We show how to do it by first year calculus. The following section in which we prove the inequality
in general is in fact a straightforward generalization of this approach.

We take the case $j=2$ with $k_1$ and $k_2$ distinct. We see from the definition that
$$F_{n,k_1,k_2}(\lambda) = (1-\lambda)^2 (n+k_1-\lambda)(n+k_2-\lambda)-k_1k_2.$$
Since $F$ is an upward facing quartic, we see that since $F(n)$ is positive and $F(n+k_1)$ is
negative, $F$ has at least two roots greater than $n$. Since $F(1)$ is negative, we see that there
are only two roots larger than $n$ and these are the roots we must sum. Let us denote them
in decreasing order, $s_1$ and $s_2$.

Now consider the quadratics
$$Q_1(\lambda)=(n+k_1-\lambda)(n+k_2-\lambda) -{k_1 k_2 \over (1-s_1)^2}.$$
and
$$Q_2(\lambda)=(n+k_1-\lambda)(n+k_2-\lambda) -{k_1 k_2 \over (1-s_2)^2}.$$
We see immediately that $Q_1(s_1)=Q_2(s_2)=0$. Moreover, it is not hard to see that
$s_1$ is the largest root of $Q_1$ and $s_2$ is the smallest root of $Q_2$.
Let $r_1$ be the largest root of $Q_2$ and $r_2$ be smallest root of $Q_1$.
 This follows from the fact
that $Q_1(n+k_1)$ and $Q_2(n+k_1)$ are both negative and $s_1$ is larger than $n+k_1$
while $s_2$ is smaller.

Now $Q_1$ and $Q_2$ differ by a constant. Since $s_1 > s_2$, we have that $Q_1$ is always
larger than $Q_2$. There for $Q_1 \triangleright Q_2$. Thus $s_1 + s_2 < r_1 + s_2$.  But
the right hand side is the sum of the roots of $Q_2$ which we may read off from the formula.
$$r_1+s_2=2n + k_1 + k_2.$$
This is the desired inequality.

 \head \S 4 The heart of the matter \endhead

As before, we let $k_1 \geq k_2 \dots \geq k_j \geq 1$ and $n > j \geq 2$.
We define
$$G_{n,k_1, \dots k_j}(\lambda)=1 + \sum_{l=1}^j {k_l \over
(n+k_l-\lambda)(1-\lambda)-k_l}.$$
Let
$$F_{n,k_1,\dots,k_j}(\lambda)=G_{n,k_1, \dots k_j}(\lambda) \prod_{l=1}^{j}
[(n+k_l-\lambda)(1-\lambda)-k_l].$$
The function $F_{n,k_1,\dots,k_j}(\lambda)$ is a polynomial of degree $2j$.
Note that $F_{n,k_1,\dots, k_j}(\lambda)$ is the determinant of the matrix
$$M_{n,k_1,\dots,k_l}=
\pmatrix  (n+k_1-\lambda)(1-\lambda) & 1 & \dots & 1 \\
          1 & (n+k_2-\lambda)(1-\lambda) & \dots & 1 \\
          \vdots  & \vdots & \vdots & \vdots         \\
          1 & 1  & \dots & (n+k_j-\lambda)(1-\lambda) \endpmatrix.$$

The goal of this section is to prove:

\proclaim{Main Lemma} The roots of the polynomial $F_{n,k_1,\dots,k_j}(\lambda)$ are
positive and real. Let us denote them by 
$s_1 \geq s_2 \geq \dots s_j \geq s_{j+1} \geq \dots \geq s_{2j}$.
Then
$$\sum_{l=1}^j s_l \leq jn+ \sum_{l=1}^j k_l.$$
\endproclaim

The fact that the roots are positive and real follows from the positivity
of the matrix $M_{n,k_1,\dots,k_l}$. However, we may also demonstrate it using
the intermediate value theorem, in a way that localizes the roots and counts
their multiplicity.

We define
$$r_l^+={n+k_l+1 + \sqrt{(n+k_l+1)^2 -4n}  \over 2},$$
and
$$r_l^{-}={n+k_l+1 -\sqrt{(n+k_l+1)^2 -4n}  \over 2},$$
to be the roots of the quadratic $(n+k_l-\lambda)(1-\lambda)-k_l$.
Notice that we have
$$0<r_1^{-} \leq r_2^{-} \leq \dots r_j^{-} < 1 < n+k_j \leq r_j^+ \leq
r_{j-1}^+ \leq \dots \leq r_1^+ <n+k_1+1.$$

It is evident that if $k_l$ appears in the list $k_1,\dots,k_j$ with
multiplicity $m$ then $r_l^+$ and $r_l^-$ appear amongst the roots
of $F_{n,k_1,\dots,k_j}$ with multiplicity $m-1$. All remaining
roots of $F_{n,k_1,\dots,k_j}$ are zeroes of $G_{n,k_1,\dots,k_j}$.

We note that $G_{n,k_1,\dots,k_j}$ has singularities at each of the
$r_l^{-}$'s and each of the $r_l^{+}$'s. The function $G_{n,k_1,\dots,k_j}(\lambda)$
approaches $+\infty$ from the left and $-\infty$ from the right at each $r_l^{-}$
and approaches $-\infty$ from the right and $+\infty$ from the left at each
$r_l^{+}$. (The preceding is because quadratics with distinct roots change signs at
each of their roots.) Thus, by the intermediate value theorem, the function
$G_{n,k_1,\dots,k_j}(\lambda)$ has a zero between $r_{l}^+$ and $r_{l-1}^+$ for
every $l>1$ for which $k_l$ and $k_{l-1}$ are distinct and a zero between 
$r_{l-1}^-$ and $r_l^-$ for every $l>1$ for which $k_l$ and $k_{l-1}$ are distinct.
Moreover, $G_{n,k_1,\dots,k_j}(n)=1-{j \over n}$ which is positive. Since
$G_{n,k_1,\dots,k_j}$ approaches $-\infty$ from the right at $r_j^{-}$ and
from the left at $r_j^{+}$, there must be a zero of $G_{n,k_1,\dots,k_j}$ between
$r_j^{-}$ and $n$ and another zero between $n$ and $r_j^{+}$. We have thus
accounted for all the roots of $F_{n,k_1,\dots,k_j}$ and indeed all the zeroes
of $G_{n,k_1,\dots,k_j}$. In particular, if we are interested in controlling the
sum of the $j$ largest roots of $F_{n,k_1,\dots,k_j}$, we are interested precisely
in the sum of those roots larger than $n$.

When $\lambda$ is larger than $n$, the frequently occuring
factor $(1-\lambda)$ is a large negative number whose order of maginitude varies
slowly. Therefore, we ``approximate" the factor $1-\lambda$ by a constant. In
other words,
we introduce, for any constant $a<1-n$, the functions
$$G_{n,k_1,\dots,k_j}^a(\lambda)=1 + \sum_{l=1}^j {k_l \over
(n+k_l-\lambda)a-k_l}.$$
and
$$F_{n,k_1,\dots,k_j}^a(\lambda)=G_{n,k_1, \dots k_j}^a(\lambda) \prod_{l=1}^{j}
[(n+k_l-\lambda)a-k_l].$$

We would like to be able to relate the organization of the roots of 
$F_{n,k_1,\dots,k_j}^a(\lambda)$ to that of $F_{n,k_1,\dots,k_j}(\lambda)$ so
we introduce a homotopy between them. More precisely, we define
$$G_{n,k_1,\dots,k_j}^{a,t}(\lambda)=1 + \sum_{l=1}^j {k_l \over (n+k_l-\lambda)(ta+
(1-t)(1-\lambda))-k_l}.$$
and
$$F_{n,k_1,\dots,k_j}^{a,t}(\lambda)=G_{n,k_1, \dots k_j}^{a,t}(\lambda)
\prod_{l=1}^{j} 
[(n+k_l-\lambda)(ta+(1-t)(1-\lambda))-k_l].$$

Notice that for every value of $t$ with $0 \leq t < 1$, we have that 
$F_{n,k_1,\dots,k_j}^{a,t}(\lambda)$ is a polynomial of degree $2j$ and we may
analyze
its roots as we did the roots of $F_{n,k_1,\dots,k_j}(\lambda)$.

We define $r_l^{a,t,+}$ and $r_l^{a,t,+}$ to be respectively the largest and
smallest
roots of the quadratic $(n+k_l-\lambda)( ta + (1-t)(1-\lambda))$. As before,we have
$$r_1^{a,t,-} \leq r_2^{a,t,-} \leq \dots \leq r_j^{a,t,-} \leq n \leq r_j^{a,t,+}
\dots \leq r_1^{a,t,+}.$$ Moreover $G_{n,k_1,\dots,k_j}^{a,t}(n) > 
G_{n,k_1,\dots,k_j}^{a,t}(n)>0.$ Thus we can locate all the roots of
$F_{n,k_1,\dots,k_j}^{a,t}(\lambda)$. It is evident that if $k_l$ appears with
multiplicity $m$ in the list $k_1,\dots,k_l$, then $r_l^{a,t,+}$ and $r_l^{a,t,+}$
are roots $F_{n,k_1,\dots,k_j}^{a,t}(\lambda)$ of with multiplicity $m-1$.
As long as $l>1$ and $k_l$ is distinct from $k_{l-1}$ there is a root of
$F_{n,k_1,\dots,k_j}^{a,t}(\lambda)$ between respectively, $r_l^{a,t,-}$ and
$r_{l-1}^{a,t,-}$ and between $r_{l-1}^{a,t,+}$ and $r_l^{a,t,+}$. Further,
there is a root between $r_{j}^{a,t,-}$ and $n$ and another root between
$n$ and $r_j^{a,t,+}$. We have now accounted for all the roots of  
$F_{n,k_1,\dots,k_j}^{a,t}(\lambda)$. In particular, their multiplicities do not
depend on $t$. Therefore, if we denote by $s_k(a,t)$, the $k$th largest root
of $F_{n,k_1,\dots,k_j}^{a,t}(\lambda)$, we see that $s_k(a,t)$ is real analytic in
$t$ for $0\leq t < 1$. (The roots never cross as we vary $t$ between 0 and 1.) 

We see moreover that $s_k(a,0)=s_k$, the $k$th largest root of 
$F_{n,k_1,\dots,k_j}(\lambda)$. As $t$ approaches $1$, the $j$ largest
roots remain bounded below by $n$ and interspersed among
the values $r_l^{a,t,+}$. It can be seen that
$$\lim_{t \longrightarrow 1} r_l^{a,t,1}=n+k_l-{k_l \over a}.$$

 Thus for $1 \leq l \leq j$ we have that
$$\lim_{t \longrightarrow 1} s_l(a,t)=s_l(a),$$
where $s_l(a)$ is the $l$th root of $F_{n,k_1,\dots,k_j}^{a}(\lambda)$.
Since $F_{n,k_1,\dots,k_j}^{a}(\lambda)$ is a polynomial of degree $j$, this
accounts for all its roots. When $l>j$, we have that
$$\lim_{t \longrightarrow 1} s_l(a,t)=-\infty.$$

Now we prove two sublemmas which combine into the proof of the main lemma.
The first relates the roots of certain polynomials $F_{n,k_1,\dots,k_j}^{a}$
to those of $F_{n,k_1,\dots,k_j}$ for certain values of $a$. The second is about
the roots of $F_{n,k_1,\dots,k_j}^{a}$ change as we vary $a$.

\proclaim{Sublemma 1} Let $1 \leq l \leq j$. Then $s_l(1-s_l) = s_l$. \endproclaim

\demo{Proof of Sublemma 1}. At the point $\lambda=s_l$, there is no difference
between $1-\lambda$ and $1-s_l$. Thus
$$F_{n,k_1,\dots,k_j}(s_l)=F_{n,k_1,\dots,k_j}^{a,t}(s_l)=F_{n,k_1,\dots,k_j}(s_l)=0.$$
Since the roots $s_k(a,t)$ never cross as we vary $t$, it must be that
$s_l$ is the $l$th root of $F_{n,k_1,\dots,k_j}(\lambda)$.
\qed \enddemo

\proclaim{Sublemma 2} We have that
$$\sum_{l=1}^j s_l(a)=jn +\sum_{l=1}^j k_l.$$
Further, whenever $a$ and $b$ are negative and $a>b$, we have that
$$F_{n,k_1,\dots,k_j}^{a} \triangleright F_{n,k_1, \dots, k_j}^{b}.$$
\endproclaim

\demo{Proof} 
The polynomial $F_{n,k_1,\dots,k_j}^a(\lambda)$ is the determinant of the product
$$\pmatrix a(n+k_1-\lambda) & 1 & \dots & 1\\
           1 & a(n+k_2-\lambda) & \dots & 1 \\
           \vdots & \vdots & \vdots & \vdots \\
           1 &  1 & \dots & a(n+k_j - \lambda) \endpmatrix
  \pmatrix k_1 & 0 & \dots & 0 \\
           0 & k_2 & \dots & 0 \\
           \vdots & \vdots & \vdots & \vdots \\
           0 & 0 & \dots & k_j \endpmatrix .$$
Thus the roots of the  polynomial $F_{n,k_1,\dots,k_j}^a(\lambda)$
are exactly the spectrum of the matrix
$$\pmatrix n+k_1 & {1 \over a} & \dots & {1 \over a} \\
           {1 \over a} & n+k_2 & \dots & {1 \over a} \\
           \vdots & \vdots & \vdots & \vdots \\
           {1 \over a} & {1 \over a} & \dots & n+k_j \endpmatrix.$$
The sum of the roots may be calculated by taking the trace. The majorization
follows from Corollary 1.3.
\qed \enddemo

Now we recursively combine Sublemma's 1 and 2 to prove the main lemma.
(We repeatedly apply $F_{n,k_1,\dots,k_j}^{1-s_l}(\lambda) \triangleright
F_{n,k_1,\dots,k_j}^{1-s_{l-1}}(\lambda)$ to the sum of the first $l-1$ roots.
Then we use sublemma 1 to identify a root of $F_{n,k_1,\dots,k_j}(\lambda).$)

$$\eqalign{ jn+\sum_{l=1}^j k_l &= s_j(1-s_j) + s_{j-1}(1-s_j) + \dots
s_1(1-s_j) \cr
                                &= s_j        + s_{j-1}(1-s_j) + \dots
s_1(1-s_j) \cr
                                &\geq  s_j + s_{j-1} (1-s_{j-1}) + \dots
s_1(1-s_{j-1}) \cr
                                &\geq     s_j + s_{j-1} + s_{j-2}(1-s_{j-2}) +
\dots s_1(1-s_{j-2}) \cr
                                & \dots  \cr
                                &\geq s_j + s_{j-1} + \dots + s_1}.$$

Which was to be shown.

\head \S 5 Proof of the main theorem \endhead.

We proceed to prove the Grone-Merris conjecture for 1-regular semibipartite graphs.
We prove the nondegenerate case first.

We first observe that $d_1^T=n+k_1 + \dots k_j$, that $d_l^T=n$ for $2 \leq l \leq n-1$
and that $d_n^T=j$. Finally, we observe that $d_l^T$ is decreasing and natural number 
valued.  Next we compare the eigenvalues of $L_g$. The $j$ largest, all larger than $n$
are the $j$ largest roots of $F_{n,k_1,\dots,k_j}$, which we call $s_1, \dots, s_j$. Following
close behind are $n-1$ eigenvalues of $n$ coming from the sum zero vectors on the
extra vertices. Coming behind then is $s_{j+1}$ which is smaller than $j$ as can be
verified by the fact that $G_{n,k_1,\dots,k_j}(j)$ is positive. All remaining eigenvalues are bounded
above by 1.

By the main lemma we have that
$$\sum_{l=1}^j s_j \leq \sum_{l=1}^j d_l^T.$$
This is the $j$th inequality for the theorem.
Since for every $l \leq j$ we have $s_l \geq n$, all the $l$th inequalities follow for $2 \leq l \leq j$.
The first inequality is simply a well known bound on the norm of $L_G$. The next $n-j-1$ inequalities
follow because $\lambda_l$ and $d_l^T$ both equal $n$ for $l$ in this range. The $n$th
inequality follows because $\lambda_n < j = d^T_n$. The remaining inequalities follow
since all other eigenvalues are bounded by 1. As soon some $d^T$ is  0, all later $d_T$'s are
zero and the remaining inequalities follow because the final inequality is an equality on trace grounds.
We have taken care of the nondegenerate case.

Now we must deal with the degenerate case. This we seem not to have dealt with at all
because we derived the polynomial $F_{n,k_1,\dots,k_j}$ under the assumption of the presence of
extra vertices. In fact even the matrix $M(j,k_1,\dots,k_j)$ does not arise as  the action of $L_G$
on functions on the graph constant on the sets $W_l$. Let us denote the set ofeigenvalues of $L_G$ on that space (with multiplicities) as $\Lambda(j,k_1,\dots,k_j)$. In fact
$$\Lambda(j,k_1,\dots,k_j)=\operatorname{spec}(M(j,k_1,\dots,k_j)) \backslash \{j\} =
\operatorname{roots}(F_{j,k_1,\dots,k_j}) \backslash \{j\} \cup \{0\}.$$
The last equality follows from continuity of the spectrum.

As is easily verified, $F_{j,k_1,\dots,k_j}$ has a double root at $j$ corresponding to its $j$th and $j+1$st
roots. One of these is removed to get $\Lambda(j,k_1,\dots ,k_l)$. Now the proof of the first $j$ inequalities in Grone Merris follows just as in the nondegenerate case, and the remaining inequalities are trivial since there is no further eigenvalue greater than 1. The theorem is proved.

\end